\documentclass[12pt]{amsart}

\usepackage[english]{babel}
\usepackage[T1]{fontenc}
\usepackage[utf8]{inputenc}

\theoremstyle{plain}
    \newtheorem{theorem}{Theorem}[section]
    \newtheorem{lemma}[theorem]{Lemma}
    \newtheorem{proposition}[theorem]{Proposition}
    
    \newtheorem{conjecture}[theorem]{Conjecture}

\theoremstyle{definition}

    \newtheorem{example}[theorem]{Example}

\usepackage[
    draft = false,
    unicode = true,
    colorlinks = true,
    allcolors = blue,
    hyperfootnotes = true,
    citecolor = red
]{hyperref}

\usepackage{amsfonts,amsmath,amssymb,amscd,amsthm,url}
\usepackage[inline]{enumitem}
\usepackage{graphicx,epsf,afterpage}
\usepackage{multicol}
\usepackage{array}
\usepackage{braket}
\usepackage{epigraph}
\usepackage{rotating}

\usepackage{xcolor}

\usepackage{ulem}

\usepackage{tikz}
\usetikzlibrary{positioning, calc, intersections, through, arrows}
\usetikzlibrary{graphs}
\usetikzlibrary{graphs.standard}
\usetikzlibrary{patterns}

\newcommand\norm[1]{\ensuremath{\left\lVert#1\right\rVert}}
\newcommand\abs[1]{\ensuremath{\left\lvert#1\right\rvert}}

\newcommand{\Geom}{\mathrm{Geom}}
\newcommand{\Uniform}{\mathrm{U}}

\newcommand{\Borel}{\mathcal{B}}

\renewcommand{\Pr}{\mathrm{P}}

\DeclareMathOperator{\Expect}{\mathbb{E}}

\DeclareMathOperator{\var}{var}

\DeclareMathOperator{\Img}{Im}

\newcommand{\R}{\ensuremath{\mathbb{R}}}

\newcommand{\Cx}{\ensuremath{\mathbb{C}}}

\newcommand{\N}{\ensuremath{\mathbb{N}}}

\renewcommand{\geq}{\geqslant}

\renewcommand{\leq}{\leqslant}

\newcounter{mcnt}

\newcounter{wordcnt}

\newcommand{\Lcal}{\mathcal{L}}

\def\refitem#1{\relax}

\begin{document}


\title{Quantum law of large numbers for Banach spaces}

\author{S. Dzhenzher and V. Sakbaev}




\begin{abstract}
    We consider random operators $\Omega \to \mathcal{L}(\ell_p, \ell_p)$
     for some $1 \leqslant p < \infty$.
     The law of large numbers is known in the case $p=2$ in the form of
     usual law of large numbers.
     Instead of sum of i.i.d. variables there may be considered the
     composition of random semigroups $e^{A_i t/n}$.
     constraints.
     We obtain the law of large numbers for the case $p \leqslant 2$.
\end{abstract}



\maketitle
\thispagestyle{empty}

\section{Introduction}

The interest to  random variables with values in a space of linear operators in Banach spaces and random processes with values in the Banach space of bounded linear operators arises in such problems of mathematical physics as the ambiguity of the quantization procedure for Hamiltonian systems \cite{Berezin, OSS-2016, Orlov2} and dynamics of quantum systems interacting with a quantum statistical ensemble \cite{GOSS-2022, Vo, VS14}. 

Random linear operators and the characterization of the distributions associated with them had being studied in \cite{Oseledets, Tutubalin, 
Skorokhod, New-23}. The law of large numbers (LLN) for a composition of i.i.d. random linear operators in a Banach space is formulated in terms of distribution of spectral characteristics of composition.

The analogue of the LLN is formulated for compositions of independent random linear operators in various operator topologies is studied as the probability estimation for the deviation of a random composition from its mean value in different topologies \cite{S16, S18, OSS19}. The LLN for operators composition is the non-commutative generalization of the LLN for the sum of independent random vector valued variables. Different versions of limit theorems for composition of independent identically distributed random operators are obtained in \cite{Berger, 
GOSS-2022, SSSh23}.

Most part of results on the LLN \cite{S16, S18, OSS19} are obtained for random linear operators acting in a Hilbert space. But the case of random operators acting in a Banach space is studied significantly less. Some results were obtained on the LLN for semigroups of bounded linear operators acting in a finite dimensional space \cite{GOSS-2022}.

This work is intended to reduce this gap. Random strongly continuous semigroups of bounded linear operators acting in Banach spaces of sequences $\ell_p,\ p\geq 1$, are considered. The properties of iterations of i.i.d. random semigroups are obtained. The LLN is obtained for compositions of random semigroups of linear operators acting in the spaces \(\ell _p, \ p\in [\,1,2\,]\).

Let us explain the interest to random linear operators in the space $\ell_1$. A random quantum dynamics can be considered either as a composition of random linear operator in a Hilbert space of quantum system or as a composition of random quantum channels in the set of quantum states \cite{Aar, Hol}.

The studying of random transforms of a Hilbert space of a quantum system leads to the analysis of compositions of random unitary operators acting in the Hilbert space. The theory of random transformation of a Hilbert space is developed in \cite{Kempe, OSS19}.

The interest for composition of random linear operators in a Banach space arises in studying of random quantum channels \cite{GOSS-2022, New-23, N-23}. To this aim we should consider random values in the set of completely positive linear mappings acting in a Banach space of nuclear operators. The analysis of random quantum channels preserving an orthogonal decomposition of unity lead to the studying the random linear operators in the space $\ell_1$. 

The structure of the article is the following.
In \S\ref{s:def} the main result (Theorem~\ref{t:lln-2}) is given with some auxiliary results (Lemmas~\ref{l:measur}, \ref{l:integ}, \ref{l:ind}).
In \S\ref{s:proof} proof of Theorem~\ref{t:lln-2} is given.
In \S\ref{appx:lcm} some futher development of the main result is given.
In \S\ref{appx:exam} one may find some examples on different forms of the LLN.

\section{Formal definitions and notations}
\label{s:def}

Let \((\Omega, \mathcal{F}, \Pr)\) be a probability space.

In this text we consider Banach spaces over some field $\mathbb{K}$ ($\R$ or $\Cx$).

For a Banach space $X$ we denote by $X^*$ the \emph{dual space} of linear maps \(X\to\mathbb{K}\). For \(x \in X\) and \(f \in X^*\) we idenfity \(\braket{f, x} := f(x)\).

Recall that for Banach spaces $X,Y$ the space $\Lcal(X,Y)$ is the Banach space of linear bounded operators $X \to Y$. We shortly denote \(\Lcal(X) := \Lcal(X,X)\).

Consider on $\Lcal(X,Y)$ the norm topology $\tau$ of open sets.
Let the \emph{Borel $\sigma$-algebra $\Borel(\tau)$} be the smallest $\sigma$-algebra of subsets of $\Lcal(X,Y)$ containing $\tau$.
A \emph{random operator} \(A\colon\Omega\to\Lcal(X,Y)\) is said to be
\begin{itemize}
    \item \emph{uniformly measurable} if the preimage of any Borel set from $\Borel(\tau)$ lies in $\mathcal{F}$; and

    \item \emph{weakly measurable} if for any \(x \in X\) and \(f \in Y^*\) the random variable $\braket{f, Ax}$ is measurable (with respect to usual Borel sets on $\mathbb{K}$).
\end{itemize}
It is obvious that the weak measurability follows from the uniform one.
(We will work mostly with weakly measurable operators. The difference between weak and uniform topologies can be seen e.g. in Example~\ref{ex:wot-not-norm}.)

For a weakly measurable operator \(A\colon \Omega\to\Lcal(X,Y)\) and \(x \in X\) we denote by
\[
    \norm{Ax}_Y\colon \Omega \to \R \quad\text{and}\quad \norm{A}_{\Lcal(X,Y)}\colon \Omega \to \R
\]
corresponding random variables.
Note that they may not be measurable. However, we have the following result.

Denote \(\N := \{1, 2, 3, \ldots\}\).
In this text we consider spaces $\ell_p$ of sequences with indexes in $\N$ for \(1 \leq p \leq \infty\).
Sometimes it is convenient to consider operators \(U \in \Lcal(\ell_p, \ell_q)\) as infinite matrices with entries \(u_{ij} \in \mathbb{K}\).

We shortly denote \(\norm{x}_p := \norm{x}_{\ell_p}\) and \(\norm{A}_p := \norm{A}_{\Lcal(\ell_p)}\).

\begin{lemma}[on measurability]
\label{l:measur}
    Let \(1 \leq p < \infty\) and \(1 \leq q \leq \infty\).
    Let \(A\colon\Omega\to\Lcal(\ell_p, \ell_q)\) be a weakly measurable operator.
    Then random variables $a_{ij}$ for any $i,j$, $\norm{Ax}_q$ for any $x \in \ell_p$, and $\norm{A}_{\Lcal(\ell_p, \ell_q)}$ are measurable.
\end{lemma}

\begin{proof}
    Take any \(i,j \in \N\). Then taking \(f \in \ell_q^*\) and \(x \in \ell_p\) such that $f_i=x_j=1$ and all other $f_k=x_k=0$ for \(k \not\in\{i,j\}\) we obtain by definition of the weak measurability that
    \[
        a_{ij} = \braket{f, Ax}
    \]
    is a measurable random variable.
    
    Take any $x \in \ell_p$.
    Then for any $i$
    \[
        (Ax)_i = \sum_{j=1}^\infty a_{ij}x_j
    \]
    is a measurable random variable as a converging series of random variables.
    So are
    \begin{align*}
        \norm{Ax}_q^q &= \sum_{i=1}^\infty \abs{(Ax)_i}^q
        \quad\text{for}\quad q < \infty,\quad\text{and}\\
        \norm{Ax}_q &= \sup_{i\in\N} \abs{(Ax)_i}
        \quad\text{for}\quad q = \infty.
    \end{align*}
    
    Now
    \[
        \norm{A}_{\Lcal(\ell_p, \ell_q)} = \sup_{\substack{x \in \ell_p \\ \norm{x} = 1}} \norm{Ax}_q
    \]
    is measurable since $\ell_p$ is separable.
\end{proof}

For a weakly measurable operator $A\colon\Omega\to\Lcal(X,Y)$ its \emph{expected value} $\Expect A \in \Lcal(X,Y)$ is a Pettis integral of $A$ on $\Omega$; i.e. for any $x \in X$ and $f \in Y^*$ we have $\Expect\braket{f, Ax} = \braket{f, (\Expect A)x}$.

The following lemma asserts that all operators we work with are Pettis integrable.

\begin{lemma}[on integrability]
\label{l:integ}
    Let \(1 \leq p < \infty\) and \(1 \leq q \leq \infty\).
    Let \(A\colon \Omega \to \Lcal(\ell_p, \ell_q)\) be a weakly measurable operator such that $\Expect\norm{A}_{\Lcal(\ell_p, \ell_q)} < \infty$.
    Then \(\Expect A \in \Lcal(\ell_p, \ell_q)\) exists and is unique, and we have
    \[
        \norm{\Expect A}_{\Lcal(\ell_p, \ell_q)} \leq \Expect \norm{A}_{\Lcal(\ell_p, \ell_q)}.
    \]
    Moreover, if we consider $A$ as an infinite matrix with random entries $a_{ij}$, then
    \[
        (\Expect A)_{ij} = \Expect a_{ij}.
    \]
\end{lemma}

\begin{proof}
    The uniqueness is proved as follows. Suppose that there are two different Pettis integrals $E_1$ and $E_2$.
    Then $(E_1 - E_2)x \neq 0$ for some $x \in \ell_p$.
    Then there is $f \in \ell_p^*$ such that \(\braket{f, (E_1 - E_2)x} \neq 0\).
    Now \(\braket{f, (E_1 - E_2)x} = \Expect\braket{f, Ax} - \Expect\braket{f, Ax} = 0\) leads to the contradiction.

    We prove the existence just for the `moreover' part.
    By Lemma~\ref{l:measur} (on measurability) for any $i,j$
    \[
        \Expect \abs{a_{ij}} \leq \Expect \norm{A}_\Lcal(\ell_p, \ell_q) < \infty,
    \]
    so we may define the operator $E \in \Lcal(\ell_p, \ell_q)$ by
    \[
        E_{ij} := \Expect a_{ij}.
    \]
    First let us show that $E$ really maps $\ell_p$ to $\ell_q$.
    Indeed, for any $x \in \ell_p$ by Jensen's inequality for $q < \infty$ we have
    \[
        \norm{Ex}^q_q =
        \sum_{i=1}^\infty \abs{\sum_{j=1}^\infty \Expect a_{ij}x_j}^q
        \leq
        \Expect \sum_{i=1}^\infty \abs{\sum_{j=1}^\infty a_{ij}x_j}^q =
        \Expect \norm{Ax}^q_q < \infty,
    \]
    and for $q = \infty$ we have
    \[
        \norm{Ex}_q =
        \sup_{i\in\N} \abs{\sum_{j=1}^\infty \Expect a_{ij}x_j} \leq
        \Expect \sup_{i\in\N} \abs{\sum_{j=1}^\infty a_{ij}x_j} =
        \Expect \norm{Ax}_q < \infty.
    \]
    This automatically gives us
    \[
        \norm{E}_{\Lcal(\ell_p, \ell_q)} \leq \Expect\norm{A}_{\Lcal(\ell_p, \ell_q)}.
    \]
    
    Now let us show that $E = \Expect A$.
    Indeed, for any $f \in\ell_q^*$ and $x \in \ell_p$
    \[
        \Expect \braket{f, Ax} = \Expect \sum_{i,j=1}^\infty a_{ij}f_ix_j = 
        \sum_{i,j=1}^\infty (\Expect a_{ij})f_ix_j =
        \braket{f, Ex},
    \]
    thus we really have defined Pettis integral $\Expect A$.
\end{proof}

Weakly measurable operators $A_1, A_2, \ldots\colon\Omega\to \Lcal(X,Y)$ are said to be \emph{independent} if for any $f_1, f_2,\ldots \in Y^*$ and $x_1, x_2, \ldots\in X$ the random variables $\braket{f_1, A_1x_1}, \braket{f_2, A_2x_2}, \ldots$ are independent.

Here and below we omit the sign $\circ$ in the composition of operators.

\begin{lemma}[on independence]
\label{l:ind}
    Let \(1 \leq p < \infty\) and \(1 \leq q \leq \infty\).
    Let \(A_1, \ldots, A_n \colon\Omega\to\Lcal(\ell_p, \ell_q)\) be independent weakly measurable operators with integrable norms.
    Then
    \[
        \Expect (A_1\ldots A_n)  = (\Expect A_1) \ldots (\Expect A_n).
    \]
\end{lemma}

\begin{proof}
    It is sufficient to prove the lemma in case $n=2$.
    Consider independent operators $A,B$ as infinite matrices with random entries $a_{ij}$ and $b_{ij}$.
    Then by definition for any $i,j,k,l$ the entries $a_{ij}$ and $b_{kl}$ are independent.
    Hence
    \[
        \Expect (AB)_{ij} = \Expect \sum_{k=1}^\infty a_{ik}b_{kj} =
        \sum_{k=1}^\infty (\Expect a_{ik}) (\Expect b_{kj}),
    \]
    and the lemma follows by the `moreover' part of Lemma~\ref{l:integ} (on integrability).
\end{proof}

We say that a weakly measurable operator $A\colon\Omega\to \Lcal(X)$ is the \emph{generator} of a \emph{semigroup} $e^{At}$ for $t \geq 0$.
A sequence $W_n(t)$ of compositions of random semigroups generated by weakly measurable operators $\Omega\to\Lcal(\ell_p)$ \emph{converges in probability in SOT of $\ell_q$ to a random semigroup $W(t)$ uniformly for $t > 0$ in any segment} if for any $x \in \ell_p$, $T > 0$ and $\varepsilon > 0$
\[
    \lim\limits_{n\to\infty} \Pr \Set{\sup\limits_{t \in [\,0,T\,]} \norm{({W_n(t) - W(t)})x}_q > \varepsilon} = 0.
\]

The following theorem is know in the case of unitary semigroups generated by self-adjoint operators. We consider the case of bounded generators.

\begin{theorem}[proved in \S\ref{s:proof}]\label{t:lln-2}
    Let $1 \leq p \leq 2$.
    Let $A\colon\Omega\to\Lcal(\ell_p)$ be a generator of a semigroup $e^{At}$ such that $\norm{A}_p$ lies in a bounded ball.
    Let $A_1, A_2, \ldots$ be a sequence of i.i.d. generators that distributions coincide with the distribution of $A$.
    Then $e^{A_1t/n} \ldots e^{A_nt/n}$ converges in probability in SOT of $\ell_2$ to $e^{\Expect At}$ uniformly for $t$ in any segment.
\end{theorem}

It is interesting to obtain the LLN in SOT of $\ell_p$ not $\ell_2$.
This leads us to the following conjecture.

\begin{conjecture}
    Under the assumptions of Theorem~\ref{t:lln-2}, $e^{A_1t/n} \ldots e^{A_nt/n}$ converges in probability in SOT of $\ell_p$ to $e^{\Expect At}$ uniformly for $t$ in any segment.
\end{conjecture}

Some steps for proof of this conjecture are made in \S\ref{appx:lcm}, see~e.g.~Theorem~\ref{t:lln-1}.


\section{Proof of Theorem~\ref{t:lln-2}}
\label{s:proof}


For \(U \in \Lcal(X,Y)\) denote by \(U^* \in \Lcal(Y^*, X^*)\) the \emph{adjoint operator}.

Recall that $\ell_p \subset \ell_q$ for $p < q$.

Let $1 \leq p \leq 2 \leq q \leq \infty$ be such that $\frac{1}{p} + \frac{1}{q} = 1$.
For a weakly measurable operator $A\colon\Omega\to \Lcal(\ell_p)$ with bounded norm denote
\[
    \var A := \Expect \bigl((A - \Expect A)^* (A - \Expect A)\bigr) \in \Lcal(\ell_p, \ell_q).
\]
Note that in order to make the definition correct in case $p=1$ we need to restrict $A^*$ to $\ell_1$ (otherwise the operator under the outer expectation may not be weakly measurable).

The following lemma is known in the case $p=q=2$.

\begin{lemma}[Chebyshev's inequality]
\label{l:cheb-l2}
    Let \(1 \leq p \leq 2 \leq q \leq \infty\) be such that \(\frac{1}{p} + \frac{1}{q} = 1\).
    Let \(A\colon\Omega\to \Lcal(\ell_p)\) be a weakly measurable operator with bounded norm.
    Then for any $\varepsilon>0$ and $x \in \ell_p$
    \[
        \Pr\!\Set{\norm{(A - \Expect A)x}_2 > \varepsilon} \leq \frac{\norm{(\var A)x}_q \norm{x}_p}{\varepsilon^2}.
    \]
\end{lemma}
\begin{proof}
    Since $A$ has the bounded norm, so \((A - \Expect A)^* (A - \Expect A)\) does. Hence $\var A$ exists and the statement is meaningful.
    Then for $\varepsilon>0$ and $x \in \ell_p$ we have
    \begin{multline*}
        \norm{(\var A)x}_q =
        \sup\limits_{\substack{y \in \ell_q^* \\ \norm{y} = 1}} \braket{y, (\var A)x} \geqslant
        \frac{1}{\norm{x}_{\ell_q^*}}\Braket{x, (\var A)x} = \\
        \frac{1}{\norm{x}_p}\Expect \Braket{(A - \Expect A)^*(A - \Expect A)x, x} = 
        \frac{1}{\norm{x}_p}\Expect \Braket{(A - \Expect A)x, (A - \Expect A)x} = \\
        \frac{1}{\norm{x}_p}\Expect \norm{(A - \Expect A)x}^2_2 \geqslant
        \frac{1}{\norm{x}_p}\varepsilon^2 \Pr\!\Set{\norm{(A - \Expect A)x}_2 > \varepsilon}.
    \end{multline*}
\end{proof}

The proof of Theorem~\ref{t:lln-2} is analogous to the proof of \cite[Theorem~2]{GOSS-2022}.

\begin{proof}[Proof of Theorem~\ref{t:lln-2}]
    Denote
    \[
        W_n(t) := e^{A_1t/n} \ldots e^{A_nt/n}.
    \]
    Fix $1 \leq p \leq 2$ and any $\varepsilon > 0$, $T > 0$ and $x \in \ell_p$.
    We need to prove that
    \[
        \lim\limits_{n\to\infty} \Pr \Set{\sup\limits_{t \in [\,0,T\,]} \norm{\left(W_n(t) - e^{\Expect At}\right)x}_2 > \varepsilon} = 0.
    \]
    The idea of the proof is to show the following two asymptotic equivalences in probability, which imply the theorem:
    \begin{align}
         \label{eq:lln}
         \lim\limits_{n\to\infty} \Pr \Set{\sup\limits_{t \in [\,0,T\,]} \norm{\left(W_n(t) - \Expect W_n(t)\right)x}_2 > \varepsilon} &= 0, \\
         \label{eq:chernov}
         \lim\limits_{n\to\infty} \Pr \Set{\sup\limits_{t \in [\,0,T\,]} \norm{\left(e^{\Expect At} - \Expect W_n(t)\right)x}_2 > \varepsilon} &= 0.
    \end{align}
    
    First let us prove~\eqref{eq:lln}.
    The idea of the proof is classical: to obtain the upper bound for the variance of $W_n(t)$ and then apply Chebyshev's equality.
    For $n \in \N$ and $t \geq 0$ denote
    \[
        \Delta_n(t) := e^{A_nt} - \Expect e^{At}.
    \]
    Then
    \[
        W_n(t) = \left(\Expect e^{At/n} + \Delta_1(t/n)\right) \ldots \left(\Expect e^{At/n} + \Delta_n(t/n)\right).
    \]
    For integer $0 \leq k \leq n$ and $1 \leq a_1 < \ldots < a_k \leq n$ denote
    \[
        F_{n;k;a_1, \ldots, a_k}(t) := (\Expect e^{At})^{a_1-1} \Delta_{a_1}(t) (\Expect e^{At})^{a_2-a_1-1} \ldots \Delta_{a_k}(t) (\Expect e^{At})^{n-a_k}.
    \]
    E.g.
    \[
        F_{n;0}(t/n) = \left(\Expect e^{At/n}\right)^n = \Expect W_n(t)
    \]
    is not a random, but a usual operator.
    So,
    \[
        W_n(t) = \sum_{k=0}^n\sum_{1 \leq a_1 < \ldots < a_k \leq n} F_{n;k;a_1, \ldots, a_k}(t/n).
    \]
    Now recall that in the case $p=1$ we additionally restrict all adjoint operators to $\ell_1$.
    Then
    \begin{multline*}
        \Expect W_n(t)^*W_n(t) =
        \sum_{k=0}^n\sum_{m = 0}^n\sum_{\substack{1 \leq a_1 < \ldots < a_k \leq n \\ 1 \leq b_1 < \ldots < b_m \leq n}} \Expect F_{n;k;a_1, \ldots, a_k}(t/n)^*F_{n;m;b_1, \ldots, b_m}(t/n) = \\
        \sum_{k=0}^n\sum_{1 \leq a_1 < \ldots < a_k \leq n} \Expect F_{n;k;a_1, \ldots, a_k}(t/n)^*F_{n;k;a_1, \ldots, a_k}(t/n),        
    \end{multline*}
    where the last equality follows by Lemma~\ref{l:ind} (on independence), since $\Expect \Delta_n(t) = 0$ and all $\Delta_n(t)$ are independent for different $n \in \N$.
    \begin{multline*}
        \var W_n(t) = \Expect W_n(t)^*W_n(t) - \Expect W_n(t)^* \Expect W_n(t) = \\
        \sum_{k=1}^n\sum_{1 \leq a_1 < \ldots < a_k \leq n} \Expect F_{n;k;a_1, \ldots, a_k}(t/n)^*F_{n;k;a_1, \ldots, a_k}(t/n).
    \end{multline*}
    
    Let $2 \leq q \leq \infty$ be such that $\frac{1}{p} + \frac{1}{q} = 1$.
    Now we ready to estimate $\norm{\var W_n(t)}_{\Lcal(\ell_p, \ell_q)}$.
    Let $\rho$ be the radius of the ball which bounds the image of $A$, i.e. $\norm{A}_p \leq \rho$. Then
    \[
        \norm{e^{At}}_p \leq e^{\rho t}.
    \]
    Hence by Lemma~\ref{l:integ} (on integrability)
    \begin{equation}
    \label{eq:norm-exp}
        \norm{\Expect e^{At}}_p \leq e^{\rho t}.
    \end{equation}
    From this and from \cite[Corollary~4.6.4 (mean value theorem)]{BogachevSmolyanov} we have
    \begin{equation}
    \label{eq:norm-delta}
        \norm{\Delta_n(t)}_p \leq \norm{e^{A_nt} - I}_p + \norm{\Expect e^{At} - I}_p \leq 2\rho t e^{\rho t},
    \end{equation}
    where $I\colon\ell_p \to \ell_p$ is the identity operator.
    Then by~\eqref{eq:norm-exp} and~\eqref{eq:norm-delta}
    \[
        \norm{F_{n;k;a_1, \ldots, a_k}(t)}_p \leq (2\rho t)^k e^{n\rho t}.
    \]
    In the case $p > 1$ we have the analogous estimate for $\norm{F_{n;k;a_1, \ldots, a_k}(t)^*}_q$.
    In the case $p=1$ we also have such an estimate since the restriction to $\ell_1$ do not increase the norm.
    Hence
    \[
        \norm{\var W_n(t)}_{\Lcal(\ell_p, \ell_q)} \leq
        \sum_{k=1}^n \binom{n}{k} \left(\frac{2\rho t}n\right)^{2k} e^{2\rho t} =
        \left(\left(1 + \left(\frac{2\rho t}{n}\right)^2\right)^n - 1\right)e^{2\rho t}.
    \]
    According to Taylor's formula in the mean-value form of the remainder for the first order, there is $s \in [\,0,2\rho t\,]$ such that
    \[
        \norm{\var W_n(t)}_{\Lcal(\ell_p, \ell_q)} \leq
        \left(1 + \frac{s^2}{n^2}\right)^{n-1} \frac{2s \cdot 2\rho t}{n} e^{2\rho t} \leq
        \frac{8\rho^2t^2}{n} e^{4\rho^2t^2 + 2\rho t} =: f(t, \rho, n).
    \]
    Then by Lemma~\ref{l:cheb-l2} (Chebyshev's equality) we have
    \[
        \Pr \Set{\sup\limits_{t \in [\,0,T\,]} \norm{(W_n(t) - \Expect W_n(t))x}_2 > \varepsilon} \leq
        \frac{f(T, \rho, n) \norm{x}_p^2}{\varepsilon^2} \xrightarrow[n\to\infty]{} 0,
    \]
    which is~\eqref{eq:lln}.

    It remains to prove equation~\eqref{eq:chernov}.
    Denote the~mapping $F\colon [\,0, +\infty)\to\Lcal(\ell_p)$ by the rule $F(t) := \Expect e^{At}$.
    By Chernov's theorem we have
    \[
        \lim\limits_{n\to\infty} \Pr \Set{\sup\limits_{t \in [\,0,T\,]} \norm{\left(e^{F'(0)t} - F\left(t/n\right)^n\right)x}_p > \varepsilon} = 0
    \]
    if (a) $F$ is continuous, (b) $F(0) = I$, (c) $\norm{F(t)} \leq e^{at}$ for some $a > 0$, and (d) for any $y\in\ell_1$ there is $F'(0)y := \lim\limits_{t \to +0} \dfrac{F(t)y - y}{t}$.
    Here~(a,b) are obvious,
    (c) holds for $a = \rho$,
    and~(d) holds for $F'(0) = \Expect A$.
    Now the equality $\Expect W_n(t) = F(t/n)^n$ and the inequality $\norm{y}_2 \leq \norm{y}_p$ for any $y \in \ell_p$ entail equality~\eqref{eq:chernov}.
\end{proof}


\section{Appropriate operators}
\label{appx:lcm}

Let $M \in \Lcal(\ell_2, \ell_1)$ be some bounded injective operator such that $M^* \in \Lcal(\ell_1, \ell_2)$ is also bounded and injective.
Denote by $\Lcal_M \subset \Lcal(\ell_1)$ the set of operators $U$ such that $\Img(UM) \subset \Img M$.
For $U \in \Lcal_M$ denote
\[
    \{U\}_M := M^{-1}UM \in \Lcal(\ell_2).
\]
For $C > 1$ denote by $\Lcal_M^C \subset \Lcal_M$ the set of operators $U$ such that
\[
    \norm{\{U\}_M}_2 \leq C\norm{U}_1.
\]

\begin{theorem}
\label{t:lln-1}
    Let $B\colon\Omega\to\Lcal_M^C$ be a generator of a semigroup $e^{Bt}$ such that $\norm{B}_1$ lies in a bounded ball.
    Let $B_1, B_2, \ldots$ be a sequence of i.i.d. generators that distributions coincide with the distribution of $B$.
    Then $\left\{e^{B_1t/n} \ldots e^{B_nt/n}\right\}_M$ converges in probability in SOT of $\ell_2$ to $\left\{e^{\Expect Bt}\right\}_M$ uniformly for $t$ in any segment.
\end{theorem}

\begin{proof}
    Let $\rho$ be the radius of the ball which bounds the image of $B$, i.e. $\norm{B} \leq \rho$.
    Since $B$ maps also in $\Lcal_M^C$, an operator $A:=\{B\}_M \colon \Omega\to\Lcal(\ell_2)$ maps in a bounded ball of the radius $C\rho$.
    Then we may apply Theorem~\ref{t:lln-2} for $A_i := \{B_i\}_M$ and $p=1$.
    Thus $\Expect A = \{ \Expect B\}_M$ with $\{U + V\}_M = \{U\}_M + \{V\}_M$ and $\{UV\}_M = \{U\}_M\{V\}_M$ for any $U,V \in \Lcal_M$
    concludes the proof.
\end{proof}



Below we explore the set $\Lcal_M^C$;
so for a while we forget about any randomness.

Recall that we consider an operator $U\in\Lcal(\ell_p)$ as an infinite matrix with entries $u_{nm}$ for $n,m\in\N$; i.e.
$(Ux)_m = \sum_n u_{nm}x_n$.

For such an infinite matrix $U = (u_{nm})$ and any integer $m > 0$ denote
\[
    f_m(U) := \sum_{n=1}^\infty \abs{u_{nm}}.
\]

\begin{lemma}
\label{l:series}
    If \(M \in \Lcal(\ell_2, \ell_1)\) is defined by \((Mx)_k = \frac{x_k}{k}\) for \(k \in \N\),
    $U\in\Lcal_M$ and for any integer $m>0$
    \[
        f_m(M^{-1}U) \leq d \norm{U}_1
    \]
    then $U \in \Lcal_M^C$ for $C = \frac{\pi d}{\sqrt{6}}$.
\end{lemma}


\begin{proof}
    For any $x \in \ell_2$ such that $\norm{x}_2 = 1$ we have
    \begin{multline*}
        \norm{\{U\}_Mx}_2^2 =
        \norm{M^{-1}UMx}_2^2 =
        \sum_{n=1}^\infty \left(\sum_{m=1}^\infty \frac{n}{m}u_{nm}x_m\right)^2 \leq \\
        \sum_{n=1}^\infty \left(\sum_{m=1}^\infty \left(\frac{n}{m}u_{nm}\right)^2\right) \left(\sum_{m=1}^\infty x_m^2\right) \leq
        \sum_{m=1}^\infty \frac{1}{m^2} \left(\sum_{n=1}^\infty n \abs{u_{nm}}\right)^2 = \\
        \sum_{m=1}^\infty \frac{f_m(M^{-1}U)^2}{m^2} \leq \frac{\pi^2d^2}{6} \norm{U}^2_1.
    \end{multline*}
\end{proof}

Before we describe another kind of operators in $\Lcal_M^C$, we give the following well known proposition.

\begin{proposition}
\label{p:norm-sup}
    For $U \in \Lcal(\ell_1)$ we have
    \[
        \norm{U}_1 = \sup\limits_m f_m(U).
    \]
\end{proposition}

\begin{proof}
    Denote $\sup_m f_m(U) =: \alpha$.
    
    The upper bound $\norm{U}_1 \leq \alpha$ follows since
    for any $x \in \ell_1$ such that $\norm{x} = 1$ we have
    \[
        \norm{Ux}_1 = \sum_{n=1}^\infty \abs{(Ux)_n} = \sum_{n=1}^\infty \abs{\sum_{m=1}^\infty u_{nm}x_m}
        \leq \sum_{n=1}^\infty\sum_{m=1}^\infty \abs{u_{nm}x_m} =
        \sum_{m=1}^\infty \abs{x_m} f_m(U) \leq \alpha.
    \]

    The lower bound $\norm{U}_1 \geq \alpha$ is proved as follows.
    For any $\varepsilon > 0$ take $m$ such that
    \[
        f_m(U) > \alpha - \varepsilon.
    \]
    Take $x \in \ell_1$ such that $x_m = 1$ and $x_k = 0$ for all $k\neq m$.
    Then
    \[
        \norm{U}_1 \geq \norm{Ux}_1 = \sum_{n=1}^\infty \abs{u_{nm}} > \alpha - \varepsilon.
    \]
    Hence $\norm{U}_1 \geq \alpha$.
\end{proof}

An operator $U \in \Lcal(\ell_p)$ is called a \emph{$d$-diagonal operator} if $u_{nm} = 0$ for all $\abs{n-m} > d$.
E.g. a $0$-diagonal operator is a diagonal operator, and vice versa.

\begin{lemma}
\label{l:d-diag}
    If $M \in \Lcal(\ell_2, \ell_1)$ is defined by $(Mx)_k = \frac{x_k}{k}$ for $k \in \N$ and
    $U \in \Lcal_M$ is a $d$-diagonal operator then $U \in \Lcal_M^C$ for $C = \sqrt{(2d+1)^3}$.
\end{lemma}

\begin{proof}
    Below we consider $u_{nm} = 0$ and $x_m = 0$ if $m < 1$ or $n < 1$.
    
    For any $x \in \ell_2$ such that $\norm{x}_2 = 1$ we have
    \begin{multline*}
        \norm{\{U\}_Mx}_2^2 =
        \norm{M^{-1}UMx}_2^2 =
        \sum_{n=1}^\infty \left(\sum_{m=n-d}^{n+d} \frac{n}{m}u_{nm}x_m\right)^2 \leq \\
        \sum_{n=1}^\infty \left(\sum_{m=n-d}^{n+d} \left(\frac{n}{m}u_{nm}\right)^2\right) \left(\sum_{m=n-d}^{n+d} x_m^2\right) =
        \sum_{m=1}^\infty x_m^2 \sum_{n=m-d}^{m+d} \sum_{k=n-d}^{n+d} \left(\frac{n}{k}u_{nk}\right)^2 \leq \\
        \sum_{m=1}^\infty x_m^2 (2d+1)^3 \sup_{k \in\N} \sum_{n=1}^\infty u_{nk}^2 \leq
        (2d+1)^3 \sup_{k \in N} f_k(U)^2 = (2d+1)^3 \norm{U}^2_1.
    \end{multline*}
\end{proof}

The following examples show that Lemmas~\ref{l:series} and~\ref{l:d-diag} are meaningful.

\begin{example}
    Consider the operator $U\in\Lcal(\ell_1)$ such that $u_{nm} = 2^{n-1}$.
    Then $f_m(U) = 2$ and $f_m(M^{-1}U) = 4$ for any $m$.
    Hence $\norm{U}_1 = 2$ by Proposition~\ref{p:norm-sup}.
    Then by Lemma~\ref{l:series} we obtain $U \in \Lcal_M^C$ for $C = \frac{2\pi}{\sqrt{6}}$.
    But it is obvious that $U \not \in \Lcal(\ell_2)$.
    Hence there are some operators not on $\ell_2$, which become operators on $\ell_2$ after putting on $\{\cdot\}_M$.
    
    Moreover, this is an example of an operator from $\Lcal_M^C$ which is not $d$-diagonal for any $d$.
\end{example}

\begin{example}
    Take bounded diagonal operator $U\in\Lcal(\ell_1)$ such that $a < u_{nn} < b$.
    It lies in $\Lcal_M^1$ since $\{U\}_M = U$.
    However Lemma~\ref{l:series} is not applicable, since $f_m(M^{-1}U) > ma$.

    Obviously there are plenty of analogous examples for $d$-diagonal operators.
\end{example}

\section{Examples for another operators}
\label{appx:exam}

In this section it is convenient to numerate sequences and matrices starting from zero, and to use Dirac bra--ket notation.
E.g. in this notation an operator $U \in \Lcal(\ell_p)$ may be written like
\[
    U = \sum_{n,m=0}^\infty u_{nm} \ket{n}\bra{m}.
\]

Recall that a random variable is a measurable map $\Omega\to\R$.
In this section we will write $\Expect \xi$ and $\var\xi$ for random variables just as for operators.
Recall that $\xi \sim \Geom(\frac{1}{2})$ means that
\[
    \Pr \Set{\xi = k} = 2^{-k}\quad\text{for}\quad k = 1, 2, 3, \ldots
\]

Recall that $I\colon \ell_1\to\ell_1$ is the identity operator.

\begin{example}
\label{ex:one-unit-lln}
    Let $\xi\sim \Geom(\frac{1}{2})$ and $A := \ket{0}\bra{\xi}\colon\Omega \to \Lcal_M$.
    Then $A$ does not maps to $\Lcal_M^C$ for any $C > 0$.
    However we will show that the following law of large numbers holds.
    Taking i.i.d. $A_1, A_2, \ldots$ with the same distribution as $A$, we have
    for any $x \in \ell_1$, $T > 0$ and $\varepsilon > 0$
    \[
        \lim\limits_{n\to\infty} \Pr \Set{\sup\limits_{t \in [\,0,T\,]} \norm{(e^{A_1t/n}\ldots e^{A_nt/n} - e^{\Expect At})x}_1 > \varepsilon} = 0.
    \]

    Indeed, we have
    \begin{align*}
        A &= \ket{0}\bra{\xi}, &\quad \Expect A &= \sum_{k = 1}^\infty 2^{-k} \ket{0}\bra{k}, \\
        A^2 &= 0, &\quad (\Expect A)^2 &= 0, \\
        e^{At} &= I + At, &\quad e^{\Expect At} &= I + \Expect At,
    \end{align*}
    and $A_iA_j = 0$ for any $i,j$, so
    \[
        e^{A_1t/n}\ldots e^{A_nt/n} =
        \left(I + A_1\frac{t}{n}\right)\ldots\left(I + A_n\frac{t}{n}\right) =
        I + \sum_{i = 1}^n A_i \frac{t}{n}.
    \]
    Thus, denoting $x =: \sum\limits_{k=0}^\infty x_k \ket{k}$ and $A_i = \ket{0}\bra{\xi_i}$ for i.i.d. $\xi_i$,
    \begin{multline*}
        \norm{(e^{A_1t/n}\ldots e^{A_nt/n} - e^{\Expect At})x}_1 =
        \norm{\left(\sum_{i = 1}^n A_i \frac{t}{n} - \Expect At\right)x}_1 = \\
        \abs{\sum_{i=1}^n x_{\xi_i} \frac{t}{n} - t\sum_{k=1}^\infty 2^{-k} x_k} =
        t \abs{\sum_{i=1}^n \frac{x_{\xi_i}}{n} - \Expect x_\xi}.
    \end{multline*}
    Hence by the Chebyshev's inequality for random variables
    \[
        \Pr \Set{\sup\limits_{t \in [\,0,T\,]} \norm{(e^{A_1t/n}\ldots e^{A_nt/n} - e^{\Expect At})x}_1 > \varepsilon} \leq
        \frac{T^2\var x_\xi}{n\varepsilon^2} \xrightarrow[n\to\infty]{} 0,
    \]
    where the last follows since $\var x_\xi = \sum_{k=1}^\infty x_k^2 2^{-k} - \left(\sum_{k=1}^\infty x_k 2^{-k}\right) < +\infty$ for $x \in \ell_1$.
\end{example}

\begin{example}
    Let $\xi\sim \Geom(\frac{1}{2})$ and $A := \xi\ket{0}\bra{\xi}\colon\Omega \to \Lcal_M$.
    Then $A$ does not maps to $\Lcal_M^C$ for any $C > 0$ and, moreover, its image does not lie in any bounded ball in $\Lcal(\ell_1)$.
    However we still have the law of large numbers.
    Indeed, analogously to Example~\ref{ex:one-unit-lln},
    \[
        \Pr \Set{\sup\limits_{t \in [\,0,T\,]} \norm{(e^{A_1t/n}\ldots e^{A_nt/n} - e^{\Expect At})x}_1 > \varepsilon} \leq
        \frac{T^2\var \xi x_\xi}{n\varepsilon^2} \xrightarrow[n\to\infty]{} 0.
    \]
\end{example}



\begin{example}
\label{ex:wot-not-norm}
    Let $\xi\sim \Uniform\{-1,1\}$ and $A := \sum\limits_{k=1}^\infty i\xi k\ket{k}\bra{k} \colon \Omega\to\Lcal_M$.
    Then
    \[
        e^{At} = \sum_{k=1}^\infty e^{i\xi kt}\ket{k}\bra{k}
        \qquad\text{and}\qquad
        \Expect e^{At} = \sum_{k=1}^\infty \cos(kt) \ket{k}\bra{k}.
    \]
    Thus we have the law of the large numbers in WOT, i.e. for any $\varepsilon > 0$, $T > 0$, $x \in \ell_1$, $z \in \ell_\infty$, and i.i.d. $A, A_1, A_2, \ldots$
    \begin{multline*}    
        \Pr \Set{\sup_{t \in [\,0,T\,]} \abs{\Braket{z,(e^{A_1t/n}\ldots e^{A_nt/n} -
        \Expect e^{At})x}}
        > \varepsilon}
        \stackrel{\text{for large enough $N$}}{\leq} \\
        \Pr \Set{\sup\limits_{t \in [\,0,T\,]} \abs{\sum_{k=1}^N z_k x_k \left(e^{itk \sum \xi_j/n} - 
        \left(\cos(kt/n)\right)^n
        \right)} > \varepsilon / 2} \xrightarrow[n\to\infty]{} 0,
    \end{multline*}
    since
    $\left(\cos(kt/n)\right)^n \rightrightarrows 1$ and
    $itk\sum \xi_j/n \stackrel{\Pr}{\rightrightarrows} 0$ uniformly for $k \leq N$ and $t \leq T$.

    But there is no law of large numbers in the norm topology.
    Indeed, for large enough $n$
    \[
        \sup_{t \in [\,0,T\,]}\norm{e^{A_1t/n}\ldots e^{A_nt/n} -
        \left(\Expect e^{At/n}\right)^n}_1 \geq
        \sup_{t \in [\,0,T\,]}\sup_{k} \abs{e^{itk\sum \xi_j/n} - \cos(kt/n)^n} \geq 1,
    \]
    since for any $T > 0$ there is $t \in [\,0,T\,]$ and $k \in \N$ such that $\abs{\cos(kt/n)^n} < 1/2$, but $\abs{e^{itk\sum \xi_j/n}} = 1$.
\end{example}

{ \bf Funding.} This work is supported by the Russian Science
Foundation under grant no. 19-11-00320.

\end{document}